\documentclass[letterpaper, 10 pt, conference]{ieeeconf}  
\usepackage{graphicx}      
\usepackage{bm}
\usepackage{amsmath}
\allowdisplaybreaks
\usepackage{amssymb}
\usepackage{xcolor}
\usepackage{mathtools}

\usepackage{booktabs}
\usepackage{multirow}
\usepackage{algorithm}
\usepackage{algpseudocode}

\IEEEoverridecommandlockouts                              
\title{\LARGE \bf
Feasibility-Aware Imitation Learning for Benders Decomposition*
}

\author{Bernard T. Agyeman$^{1}$, Zhe Li$^{1}$, Ilias Mitrai$^{2}$, and Prodromos Daoutidis$^{1}$
\thanks{*IM would like to acknowledge support from the McKetta Department of Chemical Engineering. PD would like to acknowledge support from NSF CBET (award number 2313289).}
\thanks{$^{1}$Bernard T. Agyeman, Zhe Li, and Prodromos Daoutidis are with the Department of Chemical Engineering and Materials Science, University of Minnesota, Minneapolis, MN, USA {\tt\small {\{agyem050, li003679, daout001\}@umn.edu}}}%
\thanks{$^{2}$Ilias Mitrai is with the McKetta Department of Chemical Engineering, The University of Texas at Austin, Austin, TX, USA {\tt\small imitrai@che.utexas.edu}}%
}

\begin{document}

\maketitle
\thispagestyle{empty}
\pagestyle{empty}

\begin{abstract}
Mixed-integer optimization problems arise in a wide range of control applications. Benders decomposition is a widely used algorithm for solving such problems by decomposing them into a mixed-integer master problem and a continuous subproblem. A key computational bottleneck is the repeated solution of increasingly complex master problems across iterations. In this paper, we propose a feasibility-aware imitation learning framework that predicts the values of the integer variables of the master problem at each iteration while accounting for feasibility with respect to constraints governing admissible integer assignments and the accumulated Benders feasibility cuts. The agent is trained using a two-stage procedure that combines behavioral cloning with a feasibility-based logit adjustment to bias predictions toward assignments that satisfy the evolving cut set. The agent is deployed within an agent-based Benders decomposition framework that combines explicit feasibility checks with a time-limited solver computation of a valid lower bound. The proposed approach retains finite convergence properties, as the lower bound is certified at each iteration. Application to a prototypical case study shows that the proposed method improves solution time relative to existing imitation learning approaches for accelerating Benders decomposition, while preserving solution accuracy.
\end{abstract}

\section{INTRODUCTION}
Mixed-integer optimization problems arise in a wide range of applications in engineering and control. A typical example is mixed-integer model predictive control \cite{mcallister2022advances, sager2006numerical},
in which discrete and continuous decisions must be made simultaneously to achieve a desired objective. Such problems are typically solved using branch-and-bound algorithms which iteratively explore and prune the search space. Despite significant progress in theory and algorithms, the online solution of mixed-integer optimization problems remains challenging, particularly in settings where such problems must be solved repeatedly under varying parameter realizations.

Several approaches have been proposed to reduce the solution time. For example, one can solve the problem using tailored heuristics \cite{takapoui2020simple}, warm-start the solver with past solutions \cite{marcucci2020warm}, and exploit the underlying structure with decomposition-based optimization algorithms \cite{menta2020learning, mitrai2022multicut}. Recently, machine learning (ML) approaches have been developed to reduce the online computational effort. ML approaches are particularly attractive in settings where the same optimization problem is solved repeatedly under varying parameter realizations, as past solutions can be leveraged to accelerate future computations. For example, ML models, such as neural networks and decision trees, can be used to approximate the solution to the optimization problem itself \cite{karg2020efficient, kumar2021industrial, mitrai2025discovering}, i.e., learn the mapping from the problem parameters to the optimal solution. Although such approaches can significantly reduce the solution time, the predicted solution is not always guaranteed to be optimal or feasible. While optimality is usually related to system performance, e.g., cost and/or control performance, feasibility is related to safety, i.e., constraint satisfaction. 

An alternative use of ML is to accelerate the solution algorithm itself \cite{mitrai2025accelerating, xavier2021learning}. In this approach, ML is used to either speed up calculations within an algorithm or substitute heuristic decisions. For example, ML surrogate models can be used to select the best solution strategy \cite{mitrai2024taking}, learn to select nodes in branch-and-bound solution methods \cite{nair2020solving}, and approximate/select cutting planes \cite{tang2020reinforcement, mitrai2024computationally}.

In this paper, we focus on Benders decomposition, which decomposes the problem into a mixed-integer master problem and a continuous subproblem that are solved iteratively. The master problem and subproblem are coordinated via Benders cuts, which convey optimality and feasibility information from the subproblem to the master problem. A key challenge in the online implementation of the algorithm is the repeated solution of the master problem. Solving the master problem to global optimality at early iterations can increase solution time without significantly improving the bounds~\cite{geoffrion1974multicommodity}. To accelerate the solution of the master problem, one approach is to manage and identify valuable Benders cuts using ML techniques~\cite{jia2021benders, lee2020accelerating}. More recently, reinforcement learning (RL) has been used to select the level of sub-optimality at each iteration~\cite{li2025learning}. However, these approaches still require solving a mixed-integer programming (MIP) problem at every iteration, which remains computationally demanding. To reduce solution time and reliance on the MIP solver, imitation learning and RL have been used to directly predict the solution of the master problem, where an agent selects the values of the integer variables at which Benders cuts are generated~\cite{agyeman2025graph, mak2023towards}. However, feasibility with respect to the master problem constraints is not explicitly enforced during training or deployment, and the resulting predictions may be infeasible. Moreover, due to the absence of optimality guarantees in imitation learning and RL approaches, the lower bound derived from the agent’s predictions is not guaranteed to be valid, and therefore, convergence of the resulting Benders decomposition algorithm cannot be guaranteed.

In this paper, we propose a feasibility-aware imitation learning approach to predict the values of the complicating variables in each iteration. First, the agent learns to imitate the behavior of an expert, i.e., a branch-and-bound solver, by predicting the optimal solution of the master problem at each iteration. In a second step, a feasibility-aware loss function is used to fine-tune parts of the neural network architecture to promote predictions that are feasible with respect to the accumulated Benders feasibility cuts. The agent is deployed within an agent-based Benders decomposition framework, in which explicit feasibility checks are combined with a time-limited solver computation of a valid lower bound. This ensures that the iterates of the integer variables are feasible with respect to the master problem and that lower-bound updates are valid, thereby ensuring finite convergence of the agent-based Benders decomposition framework. The proposed approach is applied to a prototypical case study, where it demonstrates a reduction in solution time without compromising solution accuracy.

\section{PRELIMINARIES}
\subsection{Problem Formulation}\label{sec:prob_form}
We consider mixed-integer nonlinear programming (MINLP) problems of the form:
\begin{align}\label{eq:minlp_form}
	\min_{\bm{x}, \bm{y}} \; &  f(\bm{x}, \bm{y}) \\
	\text{s.t.} \;& \bm{h}(\bm{x},\bm{y})  = \bm{0}, \\
	& \bm{g}(\bm{x}, \bm{y}) \leq \bm{0}, \\
	& \bm{K}\bm{y} - \bm{b} \leq 0, \label{eq:pure_binary_constraints}\\
    & \bm{x} \in \mathcal{X}, \bm{y} \in \{0,1\}^m
\end{align}
where \( \bm{x} \in \mathbb{R}^n \) are continuous variables and \( \bm{y} \in \{0,1\}^m \) are binary variables. The restriction to binary variables is without loss of generality; they arise in a wide range of applications whereas general integer variables can always be represented using binary encodings. The set \( \mathcal{X} \subseteq \mathbb{R}^n \) is defined as \(
\mathcal{X} \coloneqq \left\{ \bm{x} \in \mathbb{R}^n \mid \bm{E}\bm{x} \leq \bm{d} \right\},
\) where \( \bm{E} \in \mathbb{R}^{l \times n} \), \( \bm{d} \in \mathbb{R}^l \), \( \bm{K} \in \mathbb{R}^{o \times m} \),  and \( \bm{b} \in \mathbb{R}^o \). The function \( f:\mathbb{R}^n \times \mathbb{R}^m \to \mathbb{R} \), and the vector-valued functions \( \bm{h}:\mathbb{R}^n \times \mathbb{R}^m \to \mathbb{R}^p \) and \( \bm{g}:\mathbb{R}^n \times \mathbb{R}^m \to \mathbb{R}^q \), are assumed to be continuous, differentiable, and convex. The constraints involving only $\bm{y}$, given by constraint~\eqref{eq:pure_binary_constraints}, are referred to as \textit{pure binary constraints} throughout the paper.
The number of continuous and binary variables is denoted by $n$ and $m$, respectively, while $p$, $q$, $l$, and $o$ denote the number of equality constraints, inequality constraints, constraints defining $\mathcal{X}$, and pure binary constraints.

\subsection{Generalized Benders Decomposition (GBD)}\label{sec:gbd}
In this work, we focus on generalized Benders decomposition (GBD)~\cite{geoffrion1972generalized}, which extends classical Benders decomposition to mixed-integer nonlinear problems. The proposed approach is equally applicable to classical Benders decomposition. In particular, we adopt the GBD algorithm presented in~\cite{floudas1995nonlinear}, which incorporates both equality and inequality constraints. GBD iteratively refines upper and lower bounds on the optimal solution of problem~\eqref{eq:minlp_form}. The lower bound (LBD) is obtained from the master problem, while the upper bound (UBD) is obtained from the subproblem. In iteration $k$, the subproblem $\mathcal{S}(\bm{y}^k)$ corresponds to problem~\eqref{eq:minlp_form} with fixed binary variables $\bm{y}^k$, and is given by:
\begin{equation}
	\begin{aligned}\label{eq:minlp_sp}
		Z(\bm{y}^k) = \min_{\bm{x}} \; & f(\bm{x},\bm{y}^k)  \\
		\text{s.t.} \;& \bm{h}(\bm{x},\bm{y}^k)  = \bm{0}, \\
		& \bm{g}(\bm{x},\bm{y}^k) \leq \bm{0},\\
        & \bm{x} \in \mathcal{X}
	\end{aligned}
\end{equation}
If $\mathcal{S}(\bm{y}^k)$ is feasible, its objective value provides an UBD; otherwise, a feasibility subproblem is solved to generate a feasibility cut that excludes the current assignment. The feasibility subproblem $\mathcal{F}(\bm{y}^k)$ is given by:
\begin{equation*}
	\begin{aligned}
		\min_{\bm{x} \in \mathcal{X},\bm{\alpha}} \; & \sum_{i=1}^{q} \alpha_i \\
		\text{s.t.} \quad & \bm{h}(\bm{x},\bm{y}^k)  = 0, \\
		& {g}_i(\bm{x},\bm{y}^k) \leq \alpha_i, \quad i = 1,\dots,q,\\
		&  \alpha_i \geq 0, \quad i = 1,\dots,q.\\		
	\end{aligned}
	\label{eq:sp_infeas}
\end{equation*}
Under this decomposition, problem~\eqref{eq:minlp_form} can be equivalently reformulated as:
\begin{equation}
	\begin{aligned}\label{eq:minlp_form_mp}
		\min_{\bm{y} \in \{0,1\}^m} \; &  Z(\bm{y}) \\
		\text{s.t.} \;& \bm{K}\bm{y} - \bm{b} \leq 0. \\
	\end{aligned}
\end{equation}
Since $Z(\bm{y})$ is not available in closed form, it is approximated using hyperplanes, known as Benders optimality and feasibility cuts. The optimality cut $\mathcal{C}_o^{(k)}$ is given by:
\begin{equation}
	\mu_B \geq f(\bm{x}^k,\bm{y}) + \bm{\lambda}_k^\top \bm{h}(\bm{x}^k,\bm{y}) + \bm{\mu}_k^\top \bm{g}(\bm{x}^k,\bm{y}),
\end{equation}
where $\mu_B$ is an auxiliary variable that represents the objective of the master problem, $\bm{x}^k$ is an optimal solution of $\mathcal{S}(\bm{y}^k)$, and $\bm{\lambda}_k$, $\bm{\mu}_k$ are the associated optimal dual multipliers. We define the function $\phi_o^{(k)}(\bm{y})$ as the right-hand side (RHS) of the optimality cut, i.e.,
\[
\phi_o^{(k)}(\bm{y}) \coloneqq f(\bm{x}^k,\bm{y}) + \bm{\lambda}_k^\top \bm{h}(\bm{x}^k,\bm{y}) + \bm{\mu}_k^\top \bm{g}(\bm{x}^k,\bm{y}).
\]
The feasibility cuts $\mathcal{C}_f^{(k)}$ are given by:
\begin{equation}
	0 \geq \overline{\bm{\lambda}}_k^\top \bm{h}(\overline{\bm{x}}^k, \bm{y}) + \overline{\bm{\mu}}_k^\top \bm{g}(\overline{\bm{x}}^k,\bm{y}),
\end{equation}
where $\overline{\bm{x}}^k$ is an optimal solution of $\mathcal{F}(\bm{y}^k)$, and $\overline{\bm{\lambda}}_k$, $\overline{\bm{\mu}}_k$ are the associated dual multipliers. We define the function $\phi_f^{(k)}(\bm{y})$ as the RHS of the feasibility cut, i.e.,
\[
\phi_f^{(k)}(\bm{y}) \coloneqq \overline{\bm{\lambda}}_k^\top \bm{h}(\overline{\bm{x}}^k, \bm{y}) + \overline{\bm{\mu}}_k^\top \bm{g}(\overline{\bm{x}}^k,\bm{y}).
\]
The master problem $\mathcal{M}$ is given by:
\begin{equation}
	\begin{aligned}
		\min_{\mu_B,\, \bm{y} \in \{0,1\}^m} \; & \mu_B \\
		\text{s.t.} \quad 
		& \mu_B \ge \phi_o^{(k)}(\bm{y}), \quad \forall k \in K_{O}, \\
		& \phi_f^{(k)}(\bm{y}) \le 0, \quad \forall k \in K_{F}, \\
		& \bm{K}\bm{y}  - \bm{b} \leq 0,\\
	\end{aligned}
	\label{eq:master}
\end{equation}
where $K_O$ and $K_F$ denote the index sets of optimality and feasibility cuts. Since the number of cuts can be large, they are added adaptively during the solution process. Under standard convexity assumptions on the subproblem, the sequences of lower and upper bounds are nondecreasing and nonincreasing, respectively, and converge in a finite number of iterations.
\section{FEASIBILITY-AWARE IMITATION LEARNING}
\begin{figure}[t]
    \centering
    \includegraphics[trim={3.5cm 4.5cm 19.5cm 7.2cm}, clip, scale=0.65]{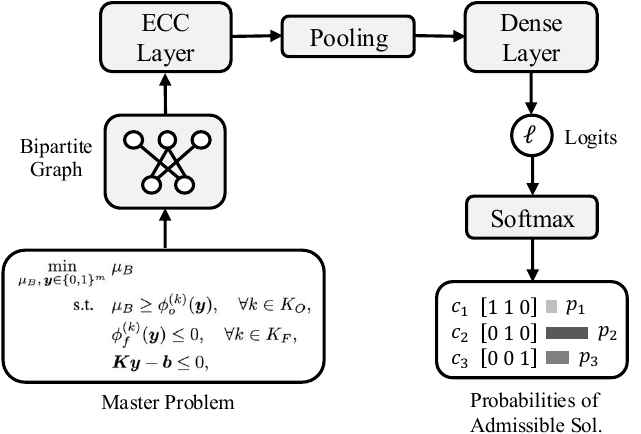}
    \caption{The architecture of the graph-based master problem agent.}
    \label{fig:diagram1}
\end{figure}

We propose a feasibility-aware imitation learning approach (see Figure~\ref{fig:diagram1}) to solve the master problem in GBD. Specifically, behavioral cloning~\cite{bain1995framework} is used to train an agent from past master problem solutions generated by an MIP solver, which serve as expert demonstrations. The agent takes as input a bipartite graph representation of the current master problem and predicts the values of $\bm{y}$.

To motivate the proposed design, we examine the structure of the master problem. By eliminating the variable $\mu_B$ in $\mathcal{M}$, problem \eqref{eq:master} can be equivalently expressed as:
\begin{equation}
	\begin{aligned}
		\min_{\bm{y} \in \{0,1\}^m} \quad & \max_{k \in K_O} \ \phi_o^{(k)}(\bm{y}) \\
		\text{s.t.} \quad 
		& \phi_f^{(k)}(\bm{y}) \le 0, \quad \forall k \in K_F, \\
		& \bm{K}\bm{y} - \bm{b} \le 0. \\
	\end{aligned}
	\label{eq:master_equiv}
\end{equation}
From problem~\eqref{eq:master_equiv}, the constraints of the master problem can be divided into two categories: (i) pure binary constraints, which are common to all instances of $\mathcal{M}$ and determine the admissible values of $\bm{y}$, and (ii) feasibility cuts, which are introduced adaptively during the GBD iterations. We leverage this structure in the design and training of the agent. Since the pure binary constraints are part of the original problem and are fixed and fully known, we enforce them by restricting the agent's action space to binary combinations that satisfy constraint~\eqref{eq:pure_binary_constraints}. In contrast, feasibility cuts are added adaptively, and the agent cannot anticipate the full set of cuts it will encounter. Therefore, their satisfaction is promoted during training through a feasibility-based logit adjustment that penalizes actions in the agent's action space that violate the current feasibility cut set. The overall framework combines a graph-based policy with a two-stage training procedure.

\subsection{Action Space Design}
The full binary action space of the master problem is given by $\mathcal{Y} \coloneqq \{0,1\}^m$, which contains \(2^m\) possible assignments. We use the pure binary constraints in constraint~\eqref{eq:pure_binary_constraints} to eliminate assignments that are infeasible \textit{a priori}. This yields the reduced action space
\[
\mathcal{Y}_{\mathrm{pb}} \coloneqq \{\bm{y}\in\mathcal{Y}\mid \bm{K}\bm{y}-\bm{b}\le 0\}.
\]
Let
\[
\mathcal{Y}_{\mathrm{pb}}=\{\bm{c}^{(1)},\ldots,\bm{c}^{(N_c)}\},
\]
where \(N_c = |\mathcal{Y}_{\mathrm{pb}}|\), \(\bm{c}^{(i)} \in \{0,1\}^m\) denotes the \(i\)-th binary assignment satisfying the pure binary constraints, and typically \(N_c \ll 2^m\). The set \(\mathcal{Y}_{\mathrm{pb}}\) defines the admissible action space of the agent, ensuring that all selected combinations \(\{\bm{c}^{(i)}\}_{i=1}^{N_c}\) satisfy the pure binary constraints.

\subsection{Agent Design}
The policy of the agent is parameterized using a graph neural network (GNN), and the learning problem is formulated as a graph-level classification task. The key components of the agent's policy are described below.

\subsubsection{Input}
The input to the agent is a bipartite graph representation of the current master problem instance, defined as
\(
\mathcal{G} \coloneqq (V, E, A_d, X_f, X_e),
\)
where $V$ and $E$ denote the sets of nodes and edges, respectively, $A_d$ is the adjacency matrix, $X_f$ is the node feature matrix, and $X_e$ is the edge feature matrix. The construction of $\mathcal{G}$ is described as follows:
\begin{enumerate}
	\item \textbf{Nodes}: The node set $V$ is partitioned into variable (binary variables) and constraint nodes (optimality and feasibility cuts). The pure binary constraints are not included as constraint nodes, as they are used to define the action space of the agent.
	
	\item \textbf{Node features}: Variable node features correspond to the values of the binary variables from the previous master problem iteration. Constraint node features include the RHS values of the constraints and a cut-type indicator, where feasibility cuts are assigned a value of 1 and optimality cuts a value of 0.
	
	\item \textbf{Edges}: An edge connects a variable node and a constraint node whenever the variable appears in that constraint.
	
	\item \textbf{Edge features}: The edge feature matrix $X_e$ contains the nonzero coefficients linking binary variables and constraints.
\end{enumerate}

\subsubsection{Graph Layers}
We employ Edge-Conditioned Convolution (ECC)~\cite{simonovsky2017dynamic}, which incorporates edge features into the message-passing process.

\subsubsection{Readout Layer}
A global sum pooling operation aggregates node embeddings into a graph-level representation.

\subsubsection{Dense Layers}
Fully connected layers are applied after pooling to increase the expressive capacity of the model.

\subsubsection{Output Layer}
The output layer is defined over the reduced action space $\mathcal{Y}_{\mathrm{pb}}$. Given a graph $\mathcal{G}$, the policy network $\pi_{\theta}$ produces a vector of logits
\[
\bm{\ell} = \pi_\theta(\mathcal{G}) \in \mathbb{R}^{N_c},
\]
where $\ell_j$ corresponds to the unnormalized score associated with the binary combination $\bm{c}^{(j)}$. The output layer therefore consists of $N_c$ units with no activation.

\subsection{Two-Stage Training Procedure}
The policy network is trained using a two-stage supervised learning approach, summarized in Algorithm~\ref{alg:feas_il}. The two-stage design allows the model to first learn to imitate expert decisions and subsequently incorporate feasibility-awareness without altering the learned graph representations.

In the first stage, given a graph $\mathcal{G}$, the policy network produces logits $\bm{\ell} \in \mathbb{R}^{N_c}$, which are converted to probabilities via the softmax function,
\[
p_j=\frac{\exp(\ell_j)}{\sum_{k=1}^{N_c}\exp(\ell_k)}, \qquad  \forall j=1,\dots,N_c.
\]
The resulting vector $\bm{p}$ defines a probability distribution over the admissible binary assignments in $\mathcal{Y}_{\mathrm{pb}}$. The policy network $\pi_{\theta}$ is trained using the cross-entropy loss
\[
\mathcal{L}_{\mathrm{CE}}^{(1)}=-\log p_{j^\star},
\]
where $j^\star$ is the index corresponding to the expert solution $\bm{y}^\star$.

The loss $\mathcal{L}_{\mathrm{CE}}^{(1)}$ encourages the agent to reproduce expert solutions obtained from past master problem solves, which are feasible and optimal with respect to the current master problem. However, the learned policy is not guaranteed to match the expert exactly at every iteration and may select alternative assignments. In the context of GBD, such deviations are not problematic: suboptimal assignments that remain feasible with respect to the master problem can still be used to generate valid cuts or serve as valid candidates for the MIP solver. Therefore, beyond matching the expert solution, it is important to ensure that the agent's actions remain feasible with respect to the feasibility cuts. This motivates the incorporation of a feasibility-based logit adjustment, which biases the policy toward selecting assignments that satisfy the current cut set. This forms the basis of the second stage.

In the second stage, building on the model learned in the first stage, the graph layers are frozen, i.e., their parameters are kept fixed, and only the parameters of the dense layers are updated. From the GBD algorithm presented in Section~\ref{sec:gbd}, the feasibility cuts in the current master problem can be written as
\[
\phi_f^{(k)}(\bm{y}) \le 0, \qquad \forall k \in K_F.
\]
For each $\bm{c}^{(j)} \in \mathcal{Y}_{\mathrm{pb}}$, an infeasibility score is computed as
\[
s_j=\sum_{k \in K_F}\max\bigl(0,\phi_f^{(k)}(\bm{c}^{(j)})\bigr), \qquad \forall j=1,\dots,N_c.
\]
Since feasibility cuts are generally affine in $\bm{y}$, the computation of the infeasibility scores reduces to evaluating linear functions over the binary combinations in $\mathcal{Y}_{\mathrm{pb}}$. Consequently, $\{s_j\}_{j=1}^{N_c}$ can be computed efficiently using matrix--vector operations.

These infeasibility scores are used to adjust the output of the neural network according to
\[
\tilde{\ell}_j=\ell_j-\omega s_j, \qquad j=1,\dots,N_c,
\]
where $\omega>0$ is a penalty parameter. The adjusted logits are converted to probabilities via
\[
\tilde{p}_j=\frac{\exp(\tilde{\ell}_j)}{\sum_{k=1}^{N_c}\exp(\tilde{\ell}_k)}.
\]
This adjustment biases the probability distribution toward combinations that satisfy the feasibility cuts. The policy network is then fine-tuned using the cross-entropy loss computed from the adjusted probability distribution
\[
\mathcal{L}_{\mathrm{CE}}^{(2)}=-\log \tilde{p}_{j^\star}.
\]

\noindent \textit{\textbf{Remark:}} The cross-entropy loss used in the first stage and the feasibility-aware loss in the second stage are complementary and could, in principle, be combined into a single training objective. However, in our experiments, the two-stage procedure was found to yield more stable training and improved performance.
\begin{algorithm}[t]
	\small
	\caption{Feasibility-Aware Imitation Learning for GBD.}
	\label{alg:feas_il}
	\begin{algorithmic}[1]
		\State \textbf{Input:} Training set $\{(\mathcal{G}_i, j_i^\star)\}$, admissible set $\mathcal{Y}_{\mathrm{pb}}$
		\State Initialize policy network $\pi_\theta$
		
		\State \textbf{Stage 1}
		\For{each minibatch}
		\State $\bm{\ell} \gets \pi_\theta(\mathcal{G})$
		\State $\bm{p} \gets \mathrm{softmax}(\bm{\ell})$
		\State $\mathcal{L} \gets -\log p_{j^\star}$
		\State Update $\pi_\theta$
		\EndFor
		\State Freeze graph  layers
		\State \textbf{Stage 2}
		\For{each minibatch}
		\State $\bm{\ell} \gets \pi_\theta(\mathcal{G})$
		\For{each $\bm{c}^{(j)} \in \mathcal{Y}_{\mathrm{pb}}$}
		\State $s_j \gets \sum_{k \in K_F} \max\!\left(0, \phi_f^{(k)}(\bm{c}^{(j)})\right)$
		\State $\tilde{\ell}_j \gets \ell_j - \omega s_j$
		\EndFor
		\State $\tilde{\bm{p}} \gets \mathrm{softmax}(\tilde{\bm{\ell}})$
		\State $\mathcal{L} \gets -\log \tilde{p}_{j^\star}$
		\State Update $\pi_\theta$
		\EndFor
	\end{algorithmic}
\end{algorithm}

\section{AGENT-BASED GBD ALGORITHM}
While the agent incorporates feasibility awareness during training, there is no guarantee that its predictions will satisfy all accumulated feasibility cuts or will be optimal with respect to the current master problem. Using the agent to directly replace the solver would therefore compromise the convergence of GBD. Infeasible predictions may lead to infeasible subproblems and repeated generation of redundant feasibility cuts, which will prevent updates to the UBD and affect finite convergence. Moreover, since the agent provides no optimality guarantees, the objective value induced by its assignment cannot be guaranteed to represent a valid lower bound on the optimal value of problem~\eqref{eq:minlp_form}, and thus cannot be used to reliably update the LBD. To address these issues, we design an approach (see Algorithm~\ref{alg:gbd_agent}) that regulates the acceptance and validation of the agent's assignment within GBD.

At each iteration, the agent predicts a binary assignment $\hat{\bm{y}}$ based on the graph representation of the current master problem. Since $\hat{\bm{y}}$ satisfies the pure binary constraints by construction, it is evaluated against the accumulated feasibility cuts by computing $\phi_f^{(k)}(\hat{\bm{y}})$ for all $k \in K_F$. If there exists $k \in K_F$ such that $\phi_f^{(k)}(\hat{\bm{y}}) > 0$, the agent's assignment is rejected, and the master problem is solved to optimality using the MIP solver.

When $\hat{\bm{y}}$ satisfies all feasibility cuts, it is accepted as the candidate solution for the current iteration. As noted, the lack of optimality guarantees means that the objective value induced by $\hat{\bm{y}}$ cannot be used to update the LBD. To obtain a valid lower bound, the MIP solver is invoked to solve the current master problem under a prescribed time budget $T_k$, with $\hat{\bm{y}}$ provided as a warm start. The primary objective of this solve is to obtain the best bound $\overline{\mathrm{LB}}$, i.e., the tightest lower bound certified by the branch-and-bound search within the time budget, which is then used to update the LBD. The time-limited solve also yields an integer solution $\overline{\bm{y}}$ with associated objective value $\overline{\mu}_B$, which is compared against the agent's assignment. In this comparison, the objective value $\hat{\mu}_B$ induced by $\hat{\bm{y}}$ under the accumulated optimality cuts is computed as
\[
\hat{\mu}_B \coloneqq \max_{k \in K_O} \phi_o^{(k)}(\hat{\bm{y}}).
\]
If $\hat{\mu}_B \le \overline{\mu}_B + \varepsilon_{\text{tol}}$, where $\varepsilon_{\text{tol}} > 0$ is a prescribed tolerance, then $\hat{\bm{y}}$ is accepted as the iterate; otherwise, $\overline{\bm{y}}$ is used, as accepting suboptimal solutions may slow convergence.

The time budget $T_k$ is scheduled based on the normalized optimality gap. Let $\mathrm{gap}_k \coloneqq \mathrm{UBD} - \mathrm{LBD}$ and $\tilde{\mathrm{gap}}_k \coloneqq \mathrm{gap}_k / \mathrm{gap}_0$. Given parameters $T_{\min}$ and $T_{\max}$, which represent the minimum and maximum allowable solve times, respectively, the time budget is selected as
\[
T_k \coloneqq T_{\min} + (T_{\max} - T_{\min})\left(1 - \tilde{\mathrm{gap}}_k\right).
\]
In early iterations, when the relaxation is weak, a small time budget is sufficient to obtain a valid lower bound. As the optimality gap decreases and the algorithm approaches convergence, a larger time budget is allocated to ensure higher accuracy.

Since the master problem is solved under a time limit, $\overline{\mathrm{LB}}$ may not improve monotonically across iterations, as in classical GBD. We therefore explicitly enforce monotonicity of the lower bound sequence as
\[
\mathrm{LBD}_k \coloneqq \max(\mathrm{LBD}_{k-1}, \overline{\mathrm{LB}}).
\]

\begin{algorithm}[t]
	\caption{Agent-based GBD algorithm.}
	\small
	\label{alg:gbd_agent}
	\begin{algorithmic}[1]
		\State \textbf{Input:} $\bm{y}^0$,  $\pi_{\theta}$,  $\mathcal{Y}_{\mathrm{pb}}$,  $\epsilon$,  $T_{\min}, T_{\max}$, $\varepsilon_{\text{tol}}$
		\State \textbf{Initialize:} $\mathrm{UBD} \gets +\infty$, $\mathrm{LBD} \gets -\infty$, $k \gets 0$, $K_F \gets \emptyset$, $K_O \gets \emptyset$
		
		\While{$\mathrm{UBD} - \mathrm{LBD} > \epsilon$}
		\State Solve subproblem $\mathcal{S}(\bm{y}^k)$
		\If{$\mathcal{S}(\bm{y}^k)$ is feasible}
		\State $\mathrm{UBD} \gets \min(\mathrm{UBD}, Z(\bm{y}^k))$
		\State Add optimality cut $\mathcal{C}^{(k)}_o$ to $\mathcal{M}$
		\State $K_O \gets K_O \cup \{k\}$
		\Else
		\State Solve feasibility subproblem $\mathcal{F}(\bm{y}^k)$
		\State Add feasibility cut $\mathcal{C}^{(k)}_f$ to $\mathcal{M}$
		\State $K_F \gets K_F \cup \{k\}$
		\EndIf
		
		\If{$\mathrm{UBD} - \mathrm{LBD} \le \epsilon$}
		\State \textbf{break}
		\EndIf
		
		\State Build graph $\mathcal{G}_k$ from the current master problem $\mathcal{M}$
		\State $\bm{\ell} \gets \pi_{\theta}(\mathcal{G}_k)$
		\State $\bm{p} \gets \mathrm{softmax}(\bm{\ell})$
		\State $j^{*} \gets \arg\max_{j} p_j$
		\State $\hat{\bm{y}} \gets \mathcal{Y}_{\mathrm{pb}}[j^{*}]$
		
		\If{$\phi_f^{(i)}(\hat{\bm{y}}) \le 0,\ \forall i \in K_F$}
		
		\State $\mathrm{gap}_k \gets \mathrm{UBD} - \mathrm{LBD}$
		
		\If{$k = 0$}
		\State $\mathrm{gap}_0 \gets \mathrm{gap}_k$
		\State $T_k \gets T_{\min}$
		\Else
		\State $T_k \gets T_{\min} + (T_{\max}-T_{\min})\left(1-\dfrac{\mathrm{gap}_k}{\mathrm{gap}_0}\right)$
		\EndIf
		
		\State Solve $\mathcal{M}$ under time limit $T_k$ with $\hat{\bm{y}}$ as warmstart
		\State Obtain $(\overline{\bm{y}}, \overline{\mu}_{B}, \overline{\mathrm{LB}})$
		
		\State Compute $\hat{\mu}_B \gets \max_{i \in K_O} \phi_o^{(i)}(\hat{\bm{y}})$
		
		\If{$\hat{\mu}_B \le \overline{\mu}_B + \varepsilon_{\text{tol}} $} 
		\State $\bm{y}^{k+1} \gets \hat{\bm{y}}$
		\Else
		\State $\bm{y}^{k+1} \gets \overline{\bm{y}}$
		\EndIf
		
		\State $\mathrm{LBD} \gets \max(\mathrm{LBD}, \overline{\mathrm{LB}})$
		
		\Else
		\State Solve $\mathcal{M}$ to obtain $(\overline{\bm{y}}, \overline{\mu}_{B}, \overline{\mathrm{LB}})$
		\State $\bm{y}^{k+1} \gets \overline{\bm{y}}$
		\State $\mathrm{LBD} \gets \max(\mathrm{LBD}, \overline{\mathrm{LB}})$
		\EndIf
		
		\State $k \gets k+1$
		\EndWhile
	\end{algorithmic}
\end{algorithm}

\textit{\textbf{Remark:}} 
Under the standard assumption that the subproblem is convex, the finite convergence of classical GBD relies on (i) the validity of the lower bound obtained from the master problem and (ii) the correct generation of feasibility and optimality cuts from the subproblem. In Algorithm~\ref{alg:gbd_agent}, feasibility of the iterates is enforced through explicit feasibility checks, ensuring that only admissible assignments for the master problem are used in the subproblem. The lower bound is computed using an MIP solver and, therefore, remains valid at each iteration. Furthermore, the generation of feasibility and optimality cuts follows the classical GBD procedure and is not modified by the agent. Consequently, {\it the proposed approach preserves the finite convergence properties of classical GBD}.

\section{CASE STUDY}
We consider a modified version of the problem from~\cite{floudas1995nonlinear}, denoted as $\mathcal{E}$:
\small
\begin{align}
	\min_{\bm{x},\bm{y}} \quad &
	\gamma_{1}y_{1} + \gamma_{2}y_{2} + \gamma_{3}y_{3} + \gamma_{4}y_{4} + \gamma_{5}y_{5}
	- 10x_{3} - 15x_{5} \nonumber\\
	& - 15x_{9} + 15x_{11} + 5x_{13} - 20x_{16}
	+ \exp(x_{3}) \nonumber\\
	& + \exp\!\left(\frac{x_{5}}{1.2}\right)
	- 60\ln(x_{11}+x_{13}+1) + 140 \tag{E1} \label{E1}\\
	\text{s.t.} \quad
	& -\ln(x_{11}+x_{13}+1) \le 0 \tag{E2} \label{E2}\\
	& -x_{3}-x_{5}-2x_{9}+x_{11}+2x_{16} \le 0 \tag{E3} \label{E3}\\
	& -x_{3}-x_{5}-0.75x_{9}+x_{11}+2x_{16} \le 0 \tag{E4} \label{E4}\\
	& x_{9}-x_{16} \le 0 \tag{E5} \label{E5}\\
	& 2x_{9}-x_{11}-2x_{16} \le 0 \tag{E6} \label{E6}\\
	& -0.5x_{11}+x_{13} \le 0 \tag{E7} \label{E7}\\
	& 0.2x_{11}-x_{13} \le 0 \tag{E8} \label{E8}\\
	& \exp(x_{3}) - U y_{1} \le \rho_{1} \tag{E9} \label{E9}\\
	& \exp\!\left(\frac{x_{5}}{1.2}\right) - U y_{2} \le \rho_{2} \tag{E10} \label{E10}\\
	& 1.25x_{9} - U y_{3} \le 0 \tag{E11} \label{E11}\\
	& x_{11}+x_{13} - U y_{4} \le 0 \tag{E12} \label{E12}\\
	& -2x_{9}+2x_{16} - U y_{5} \le 0 \tag{E13} \label{E13}\\
	& y_{1}+y_{2}=1 \tag{E14} \label{E14}\\
	& y_{4}+y_{5}\le 1 \tag{E15} \label{E15}\\
	& \bm{y}^{T}=(y_{1},y_{2},y_{3},y_{4},y_{5}) \in \{0,1\}^{5}, \nonumber\\
	& \bm{x}^{T}=(x_{3},x_{5},x_{9},x_{11},x_{13},x_{16}) \in \mathbb{R}^{6}, \nonumber\\
	& \bm{a}^{T}=(0,0,0,0,0,0), \quad
	\bm{b}^{T}=(2,2,2,\text{--},\text{--},3), \nonumber\\
	& \bm{a} \le \bm{x} \le \bm{b}, \nonumber\\
	& U \in [6,14], \quad \rho_{1},\rho_{2} \in [0,2], \nonumber\\
	& \gamma_{1},\gamma_{2},\gamma_{3},\gamma_{4} \in [1,39], \quad \gamma_{5} \in [1,7]. \nonumber
\end{align}
\normalsize
\subsection{Data Generation}
Expert data was obtained by sampling 300 distinct combinations of the parameters $\gamma_{1}, \gamma_{2}, \gamma_{3}, \gamma_{4}, \gamma_{5}, U, \rho_{1}, \text{ and } \rho_{2}$ within their respective ranges. For each parameter set, the resulting MINLP was solved using GBD, where Gurobi 12.0.1~\cite{gurobi} was used to solve the master problem and Ipopt 3.12.6~\cite{wachter2006ipopt} was used to solve the subproblem. The final dataset consisted of 24,638 data points.

\subsection{GNN Architecture and Training}
The initial action space of the master problem of $\mathcal{E}$ consists of $2^5 = 32$ binary combinations. By enforcing the pure binary constraints in \eqref{E14} and \eqref{E15}, this set was reduced to 12 feasible combinations, which defined the agent's action space. The policy network consisted of three ECC layers, each with 64 hidden units and ReLU activations, followed by two fully connected layers of sizes 64 and 32 with ReLU activations. The output layer consisted of 12 units corresponding to the feasible combinations. In the first stage, the model was trained using the Adam optimizer with a learning rate of $10^{-3}$, a batch size of 8, and 20 epochs. In the second stage, the  dense layers were fine-tuned using a learning rate of $10^{-4}$, a batch size of 8, 20 epochs, and a penalty parameter of $\omega = 0.1$. The training set consisted of 23,409 data points.

\subsection{Agent Evaluation}
The agent was evaluated on 1,232 data points, of which 877 instances contained feasibility cuts. The evaluation assessed the agent's ability to predict assignments that satisfied all feasibility cuts associated with each data instance in the test set, as well as its ability to match expert solutions. To assess the benefits of the two-stage training procedure, the proposed agent was compared with a baseline corresponding to the model obtained after the first training stage. The evaluation also includes the graph-based imitation learning approach  developed in~\cite{agyeman2025graph}, where an agent is trained to predict the values of the binary variables based on a bipartite graph representation of the master problem. In that approach, training is based solely on behavioral cloning, and the binary variables are predicted independently without explicitly accounting for feasibility with respect to the pure binary constraints or the accumulated Benders feasibility cuts. The agent resulting from that approach is termed independent IL in this paper. Its architecture is similar to that of the agent in the proposed approach, with the exception of the output layer, which consists of 5 units with sigmoid activation functions to reflect the 5 binary variables in problem $\mathcal{E}$. The independent IL agent was trained using the same hyperparameters as the proposed approach. Following~\cite{agyeman2025graph}, the outputs of the independent IL agent were converted to binary values using thresholds of 0.75 and 0.25.

Furthermore, 30 new distinct combinations of parameters $\gamma_{1}, \gamma_{2}, \gamma_{3}, \gamma_{4}, \gamma_{5}, U, \rho_{1}, \text{ and } \rho_{2}$ were generated. For each parameter set, the resulting problem was solved using the agent-based GBD according to Algorithm~\ref{alg:gbd_agent}. The same problem instances were also solved using the agent resulting from the independent IL approach under the same deployment scheme of Algorithm~\ref{alg:gbd_agent}. Classical GBD was applied to the same problem instances, and the approaches were compared in terms of the average time required to solve the master problem, the average number of iterations, the average subproblem solution time, and the agreement between the solutions returned by each approach and the true solution, obtained by solving the MINLP using the BONMIN 1.8.9~\cite{bonami2008bonmin} solver.  For the agent-based approaches (proposed and independent IL), we assessed the frequency with which the agent's action was accepted as the iterate and the frequency with which the time-limited master problem solved produced a better solution. For this evaluation, $T_{\min}$, $T_{\max}$, $\varepsilon$, and $\varepsilon_{\text{tol}}$ were set to 0.1 s, 0.5 s, 0.001, and $10^{-6}$, respectively. All simulations were performed on a machine equipped with an Intel Core i7-14700 processor (20 cores), 32 GB RAM, running Ubuntu 24.04.
\subsection{Results and Discussion}
\begin{table}[t]
	\centering
	\caption{Comparison of baseline, independent IL and proposed agent.}
	\label{tab:il_results}
	\begin{tabular}{lcc}
		\toprule
		Method & Exact Match & Feasibility \\
		\midrule
		Baseline & 847/1232 (68.75\%) & 781/877 (89.05\%) \\
		\midrule
		Independent IL~\cite{agyeman2025graph} & {257/1232 (20.86\%)} & {175/877 (19.95\%)} \\
        \midrule
		Proposed & \textbf{917/1232 (74.43\%)} & \textbf{822/877 (93.73\%)} \\
		\bottomrule
	\end{tabular}
\end{table}
\begin{table*}[t]
	\centering
	\caption{Comparison between classical GBD and learning-augmented variants.}
	\label{tab:gbd_results}
	\begin{tabular}{lccc}
		\toprule
		Metric & Classical GBD & Independent IL~\cite{agyeman2025graph} & Proposed \\
		\midrule
		Avg. iterations & 8.83 $\pm$ 1.21 & 8.83 $\pm$ 1.21 & 8.83 $\pm$ 1.21 \\
		Avg. MP$^{\star}$ time (s) & 0.400 $\pm$ 0.065 & 0.307 $\pm$ 0.075 & \textbf{0.139 $\pm$ 0.057} \\
		Avg. SP$^{\ast}$ time (s) & 1.240 $\pm$ 0.131 & 1.218 $\pm$ 0.276 & 1.238 $\pm$ 0.318 \\
		\midrule
		Agent-selected iterate$^{\dagger}$ & -- & 2.07 / 8.83 (23\%) & 6.27 / 8.83 (71\%) \\
		Solver-selected iterate$^{\ddagger}$ & -- & 0.20 / 8.83 (2\%) & 2.47 / 8.83 (28\%)\\
		MP$^{\star}$ time reduction & -- & 23.3\% & \textbf{65.3\%} \\
		Solution match & 30/30 (100\%) & 30/30 (100\%) & 30/30 (100\%) \\
		\bottomrule
	\end{tabular}\\
	\vspace{3pt}
	\footnotesize{\textit{Note:} 
	$^{\star}$ MP: master problem; 
	$^{\ast}$ SP: subproblem; 
	$^{\dagger}$ Agent's assignment is used as the iterate; 
	$^{\ddagger}$ Solver's solution replaces the agent's assignment. 
	All time and iteration metrics are reported as mean $\pm$ standard deviation.}
\end{table*}
Table~\ref{tab:il_results} presents a comparison between the proposed approach, the baseline, and the independent IL approach. The proposed approach improves the ability of the agent to propose assignments that satisfy feasibility cuts: out of the 877 instances with feasibility cuts, the proposed approach satisfies 822 instances (93.73\%), compared to 781 instances (89.05\%) for the baseline and 175 instances (19.95\%) for the independent IL approach. In terms of matching the expert solutions, the proposed policy also outperforms both alternatives, achieving an exact match rate of 74.43\%, compared to 68.75\% for the baseline and 20.86\% for the independent IL approach.

The improvement in feasibility satisfaction of the proposed approach over the baseline highlights the effectiveness of the feasibility-aware loss introduced in the second stage of the training process. Furthermore, the improvement in exact match indicates that this loss, while encouraging feasibility with respect to the feasibility cuts, also enhances the ability of the agent to produce actions that match the expert solutions. 
The comparison between the proposed approach and the independent IL approach highlights the importance of incorporating the structure of the master problem into the learning process. In the independent IL approach, feasibility with respect to the constraints of the master problem is learned implicitly from data. In contrast, in the proposed approach, the learning problem is restricted to binary combinations that satisfy the admissible assignment constraints, and feasibility with respect to the accumulated cuts is incorporated during training via the feasibility-based logit adjustment. This reduces the complexity of the learning problem and leads to significantly improved feasibility satisfaction and overall predictive performance. 
   
Table~\ref{tab:gbd_results} presents a comparison between classical GBD and the agent-based approaches. All approaches require the same number of iterations on average and recover the true solution in all 30 instances. In terms of computational performance, the independent IL approach achieves a reduction in the average master problem solution time, from 0.400 s to 0.307 s (23.3\%), while the proposed approach achieves a significantly larger reduction to 0.139 s (65.3\%). On average, the proposed agent's action is selected as the iterate in 71\% of the iterations, compared to 23\% for the independent IL approach.  Furthermore, in 28\% of the iterations, the time-limited master problem solve for valid lower-bound computation yields a solution that improves upon the agent's assignment and is therefore used to update the iterate; this occurs in only 2\% of the iterations for the independent IL approach. Overall, the proposed agent produces feasible assignments in 99\% of the iterations, compared to 25\% for the independent IL approach. The average subproblem solution time remains unchanged across all approaches. Importantly, whenever the proposed agent's assignment is accepted, it coincides with the solution obtained by solving the master problem to optimality with Gurobi.

\section{CONCLUSIONS}
This paper proposes a feasibility-aware imitation learning framework for accelerating the master problem in Benders decomposition. We show that existing imitation learning approaches for master problem acceleration may produce assignments that violate the constraints of the master problem. To address this, the proposed method leverages both the constraints that define admissible assignments and the accumulated feasibility cuts of the master problem. In particular, the agent’s action space is restricted to admissible assignments, while a two-stage training procedure introduces a feasibility-aware adjustment that encourages satisfaction of the feasibility cuts. The trained agent is integrated within an agent-based GBD framework, where its predictions are screened through feasibility checks and used to guide the selection of candidate iterates. A time-limited solver is employed to compute valid lower bounds, ensuring that the algorithm retains the finite convergence properties of classical GBD. Numerical experiments demonstrate that the proposed method improves both feasibility satisfaction and agreement with expert solutions compared to existing imitation learning approaches for GBD acceleration, while reducing the computational effort required to solve the master problem.

\bibliographystyle{IEEEtran}
\bibliography{refs}
\end{document}